\documentclass[graybox]{svmult}

\usepackage{type1cm}
\usepackage{makeidx}
\usepackage{graphicx}
\usepackage{multicol}
\usepackage[bottom]{footmisc}
\usepackage[colorlinks=true, linkcolor= black, linktocpage=true]{hyperref}
\usepackage{newtxtext}
\usepackage[varvw]{newtxmath}

\makeindex

\begin{document}

\title*{Linear Landau damping for a two-species Vlasov-Poisson system for electrons and ions}
\author{Lena Baumann, Marlies Pirner}
\institute{Lena Baumann \at Department of Mathematics, University of Würzburg, Emil-Fischer-Str. 40, 97074 Würzburg, Germany \\ \email{lena.baumann@uni-wuerzburg.de} \and Marlies Pirner \at Angewandte Mathematik Münster: Institut für Analysis und Numerik, 
Fachbereich Mathematik und Informatik der Universität Münster, 
Orléans-Ring 10, 
48149 Münster, Germany \\ \email{marlies.pirner@uni-muenster.de}
}
\maketitle

\abstract{This paper concerns the linear Landau damping for the two species Vlasov-Poisson system for ions and electrons near Penrose stable equilibria. The  result is an extension of the result on the one species Vlasov-Poisson equation by  Mouhout and Villani \cite{Villani, Villani_steps}. Different from \cite{Villani, Villani_steps}, we do not describe the ions as a background species but as a species which is also described by a separate Vlasov equation. We show an exponential decay of the electric energy for the linearised system near Penrose stable equilibria.}

\section{Introduction}
\label{sec:1}

We consider a plasma consisting of electrons denoted by the index $e$ and one species of ions denoted by the index $i$. The state of the two species mixture will be described by two distribution functions $f_e(x,v,t), f_i(x,v,t)\geq 0,$ $~x \in \mathbb{T}, ~ v \in \mathbb{R}, t>0$. We describe the time evolution by a system of Vlasov equations
\begin{align}
\begin{split}
&\partial_t f_e + v \partial_x f_e - \frac{e}{m_e} E[f_e, f_i] \partial_v f_e = 0,\\
&\partial_t f_i + v \partial_x f_i + \frac{e}{m_i} E[f_e, f_i] \partial_v f_i = 0,
\end{split}
\label{equ}
\end{align}
with corresponding initial data. Here, $e$ denotes the elementary charge and $m_i, m_e$ the mass of the corresponding particles. As the right-hand sides of (\ref{equ}) are equal to zero this set of equations describes a collision-less regime. The charge density $e \rho(x,t)$ is given by
\begin{align}
e ~ \rho(x,t)= e ~ \int_{- \infty}^{\infty} ( f_i(x,v,t) - f_e(x,v,t))\, \mathrm{d}v.
\label{chargeden2}
\end{align}
We assume that the magnetic field is neglectable compared to the electric field and that the electric field $E(x,t)$ can be written as the gradient of a potential $\phi(x,t)$, i.e. $E(x,t)=-\nabla \phi (x,t)$. Together with the Maxwell equation $\partial_x E(x,t) = \rho(x,t)$ we get that $\phi$ solves the Poisson equation $-\Delta \phi (x,t)= \rho (x,t)$. Under these assumptions we talk about the Vlasov-Poisson system.\\

\noindent Following \cite{Villani}, we assume that $E(x,t)$ is of the form 
\begin{align}
\begin{split}
E[f_e,f_i](x,t)&= - \int \partial_y W(y)  e ~ \rho(x-y,t) \, \mathrm{d}y \\ &= - \int \partial_y W(y)~ e \int_{-\infty}^{\infty} (f_i(x-y,v,t) - f_e(x-y,v,t))\, \mathrm{d}v \mathrm{d}y \\ &= - \int \int_{-\infty}^{\infty} \partial_y W(y)~ e ~(f_i(x-y,v,t) - f_e(x-y,v,t))\, \mathrm{d}v \mathrm{d}y,
\end{split}
\label{elec}
\end{align}
where $W$ is a Green's function to the Poisson equation on $\mathbb{T}$. Here we assume that the plasma is globally neutral $(\int \rho dx=0).$ More information on the Poisson equation in the periodic setting can be found in section $2$ of \cite{Villani}. \\ 

\noindent In this paper, we are concerned with the phenomenon of Landau damping. It is a result on the qualitative behaviour of the Vlasov-Poisson equation. In an electronic plasma a small initial disturbance from the equilibrium distribution gives rise to an oscillating electric field. In the collision-less limit Landau discovered that this self-consistent electric field is actually damped. This damping property of the Vlasov equation is referred to as Landau damping \cite{Landau}. Landau derived this phenomenon without from nowadays perspective giving a profound mathematical explanation. First experimental results approving its practical validity were given about 20 years later \cite{Malmberg}. Different physical interpretations of how this damping occurs can be found in \cite{LandauLifshitz, vanKampen, Villani_steps}. In numerical studies Landau damping serves as a practical test case to validate different numerical implementations of the Vlasov-Poisson system and is therefore still commonly used, see e.g. \cite{FilbetSonnendruecker, Filbet, Einkemmer}.\\

\noindent A mathematically more rigorous treatment of the derivation of Landau damping for the one species equation \cite{Landau} has been done by Villani in \cite{Villani, Villani_steps}. The considerations in our paper are an extension of the result presented by Villani to a mixture of ions and electrons. In this paper, we do not treat the species of ions as a background species but also describe the time evolution of the ions, additionally. 
According to the ideas of Landau one can show that when one linearises the Vlasov equations of system \eqref{equ} around a homogeneous equilibrium, it is possible to analyse the stability and the asymptotic behaviour for the linearised equations. This analysis is carried out in the following sections in a rigorous way inspired by the linearised one species case done by Villani in \cite{Villani_steps}. \\

\noindent The paper is structured as follows. We begin with the existence of equilibrium solutions in section \ref{sec6.2.2.3}. In section \ref{sec6.2.2.4} we perform the linearisation around an equilibrium distribution. In section \ref{sec6.2.2.5} we  derive analytically the property of Landau damping from the linearised equations. The main part of this work is taken from \cite{Pirner}. Section \ref{Outlook} gives a short summary and outlook.

\section{Existence of equilibria and linearisation of the Vlasov model}

\noindent In this section we want to find equilibrium solutions to \eqref{equ} and perform a linearisation around these equilibrium distributions.

\subsection{Existence of equilibria}
\label{sec6.2.2.3}
We start by giving the definition of an equilibrium solution for our two species Vlasov-Poisson system. We call a pair of functions $(f_i^{equ},f_e^{equ})$ an equilibrium solution to \eqref{equ} if and only if $(f_i^{equ},f_e^{equ})$ satisfy \eqref{equ} and $\partial_t f_i^{equ} = \partial_t f_e^{equ} =0$. Then we are able to prove that there exists an equilibrium solution to \eqref{equ}.

\begin{theorem}[Existence of equilibrium solutions]
There exists at least one equilibrium solution to \eqref{equ}.
\label{exequ}
\end{theorem}
\begin{proof}
We prove this assertion by giving an example. Any pair of distributions $(f_i(x,v), f_e(x,v))=(f_i^{equ}(v), f_e^{equ}(v))$ defines an equilibrium solution. Obviously the distributions are independent of $t$. So the only thing we have to prove is that they are a solution to \eqref{equ}. Since they are also independent of $x$, we have $v \partial_xf_i^{equ} = v \partial_x f_e^{equ} =0$ for every $v\in \mathbb{R}.$ Then the charge density $\rho^{equ}$ given by \eqref{chargeden2} associated to $f_i^{equ}$ and $f_e^{equ}$ is constant in $x$. So the corresponding forces vanish since $$\int \partial_y W(y) \rho^{equ}(x-y,t) \, \mathrm{d}y = \int W(y) \partial_x \rho^{equ}(x-y,t) \, \mathrm{d}y=0.$$
Hence, any pair of such distributions defines an equilibrium solution.
\end{proof}

\subsection{Linearisation of the Vlasov model}
\label{sec6.2.2.4}
For the linearisation of system \eqref{equ} we assume that we are near an equilibrium %and can write
\begin{align}
\begin{split}
f_e(x,v,t)= f_e^{equ}(v) + h_e(x,v,t), \\
f_i(x,v,t)=f_i^{equ}(v) + h_i(x,v,t),
\end{split}
\label{ans}
\end{align}
where $h_e,h_i$ are small deviations. %\textcolor{red}{The meaning of "small" will be specified in a moment.} 
We assume that the equilibria are independent of $x$ and they satisfy the condition of quasi-neutrality
\begin{align}
\int_{-\infty}^{\infty} f_i^{equ} (v) \, \mathrm{d}v = \int_{-\infty}^{\infty} f_e^{equ} (v) \, \mathrm{d}v,
\label{quasineu}
\end{align}
which means that in equilibrium the resulting charge $\rho$ is equal to zero. Such an equilibrium exists according to Theorem \ref{exequ}.
Inserting ansatz \eqref{ans} into system \eqref{equ} and equation \eqref{elec} and neglecting quadratic terms in $h_i$ and $h_e$, respectively, we obtain 
\begin{align}
\begin{split}
&\partial_t h_e + v \partial_x h_e - \frac{e}{m_e} E[h_e, h_i] \partial_v f_e^{equ} = 0,\\
&\partial_t h_i + v \partial_x h_i + \frac{e}{m_i} E[h_e, h_i] \partial_v f_i^{equ} = 0,
\end{split}
\label{equ2}
\end{align}
where 
\begin{align}
\begin{split}
E[h_e,h_i]&= - \int \int_{-\infty}^{\infty} \partial_y W(y) [h_i(x-y,v,t) - h_e(x-y,v,t)] \, \mathrm{d}v  \mathrm{d}y.
%\mathcolor{red}{&= - \partial_y W \ast_{x,v} h.}
\end{split}
\label{elec2}
\end{align}
We now want to show the following. For a small displacement of the electrons and the ions from equilibrium the electric fields act on the electrons and ions as a restoring force. Due to the restoring force, standing density waves are possible with a fixed frequency. So we presume to see that the densities of the deviations are oscillating in time. This oscillation is possibly damped depending on the physical configuration. For more details on the physics behind the phenomenon of Landau damping see for example \cite{Herr}. We want to show that this property is contained in our model \eqref{equ2} and find criteria when these oscillations are damped. Since the electric field depends directly via \eqref{elec2} on the densities and the electric field is related to the electric energy $W_{el}$ via $W_{el}= \frac{1}{2} E^2$, we also expect to  observe a damping when considering the electric energy $W_{el}$ or simply $E^2$.

\section{Linear Landau damping for the two species model}
\label{sec6.2.2.5}
In this section we want to derive the property of linear Landau damping for our two species model. We start with applying the method of characteristics followed by a Fourier transform of the distribution functions and the densities. Thereafter, we give concrete convergence rates of the densities and show the linear Landau damping near Landau-Penrose stable equilibria. 

\subsection{Methods of characteristics}
\label{sec6.2.2.6}
The goal of our analysis is to find criteria for damping and show the property of Landau damping in an analytical way starting from equations \eqref{equ}. Analogously to \cite{Villani} we apply the methods of characteristics to \eqref{equ2} and obtain
\begin{align}
\begin{split}
h_e(x,v,t)= h_{e,in}(x-vt,v) + \frac{e}{m_e} \int_0^t E[h_i,h_e](\tau, x-v(t-\tau)) \partial_v f_e^{equ}(v) \, \mathrm{d}\tau, \\
h_i(x,v,t)= h_{i,in}(x-vt,v) - \frac{e}{m_i} \int_0^t E[h_i,h_e](\tau, x-v(t-\tau)) \partial_v f_i^{equ}(v) \, \mathrm{d}\tau,
\end{split}
\label{char}
\end{align}
where $h_{e,in}(x,v) = h_e(x,v,0)$ and $h_{i,in}(x,v) = h_i(x,v,0)$ denote the inital values. In the next step we shall apply the Fourier transform to this set of equations.

\subsection{Fourier transform of the distribution functions}

Similar as in \cite{Villani_steps}, we take the Fourier transform  in both $x$ and $v$ of both equations \eqref{char} and do the substitutions $y=x-vt$ in the part with the initial data and $y=x-v(t- \tau)$ in the term with the time integration (source term). We use that the Fourier transform of the source terms separates into the product of the Fourier transform in $x$ of the force term $E[h_i,h_e]$ and the Fourier transform in $v$ of the derivative of the equilibrium function $\partial_v f_i^{equ}$ and $\partial_v f_e^{equ}$, respectively. Then we use for the electric field that the Fourier transform of a convolution leads to a product of two Fourier transforms. We get rid of the derivatives by doing integration by parts in $v$. All in all, we obtain
\begin{align}
\begin{split}
\widetilde{h}_e(k,\eta,t)&= \widetilde{h}_{e,in}(k,\eta+kt)\\& +\frac{e}{m_e} 4 \pi^2 \widehat{W}(k) \int_0^t \widehat{\rho}(k, \tau) \overline{f}_e^{equ}(\eta+ k(t-\tau)) k \cdot [\eta + k(t-\tau)] \, \mathrm{d}\tau, \\
\widetilde{h}_i(k,\eta,t)&= \widetilde{h}_{i,in}(k,\eta+kt)\\& - \frac{e}{m_i} 4 \pi^2 \widehat{W}(k) \int_0^t \widehat{\rho}(k, \tau)\overline{f}_i^{equ}(\eta+ k(t-\tau)) k \cdot [\eta + k(t-\tau)] \, \mathrm{d}\tau,
\end{split}
\label{fou2}
\end{align}
where $k$ and $\eta$ denote the new variables after the transformations and $\rho= \rho^{h_i} - \rho^{h_e}$, 
$\rho^{h_i}=\int h_i \, \mathrm{d}v$, $\rho^{h_e} = \int h_e \, \mathrm{d}v$. The notation $\widetilde h$ denotes the Fourier transform of a function $h(x,v,t)$ in both $x$ and $v$ variable whereas $\widehat h$ denotes the Fourier transform of $h$ only in the $x$ variable and $\overline{h}$ denotes the Fourier transform of $h$ only in the $v$ variable.\\

\noindent For $\eta=0$ we obtain the Fourier transform of the densities. It remains the integration with respect to $v$ over the function $h_i$ and $h_e$, respectively, which gives the density. We get the system
\begin{align}
\begin{split}
\widehat{\rho}^{h_e}(k,t)= \widetilde h_{e,in}(k,kt) +  \int_0^t K_e^0(t-\tau,k) [\widehat{\rho}^{h_i}- \widehat{\rho}^{h_e}](k, \tau) \, \mathrm{d}\tau, \\
\widehat{\rho}^{h_i}(k,t)= \widetilde h_{i,in}(k,kt) +  \int_0^t K_i^0(t-\tau,k) [\widehat{\rho}^{h_i}- \widehat{\rho}^{h_e}](k, \tau) \, \mathrm{d}\tau,
\end{split}
\label{den}
\end{align}
where 
$$K_e^0(t-\tau,k)=\frac{e}{m_e} 4 \pi^2 \widehat{W}(k)f_e^{equ}( kt) |k|^2 t,$$
$$K_i^0(t-\tau,k)=-\frac{e}{m_i} 4 \pi^2 \widehat{W}(k)f_i^{equ}( kt) |k|^2 t.$$
We call the functions $K_{i}^0, K_{e}^0$ the kernels of the integrals. In the next subsection we want to show the convergence rates of the densities as given in system \eqref{den}.

\subsection{Convergence rates of the densities}
Assume that the equilibrium is a global one which means $f_i^{equ}=f_e^{equ}$. It exists according to the proof of Theorem \ref{exequ}. Then 
$$K_e^0(t-\tau,k)=- \frac{m_i}{m_e} K_i^0(t-\tau,k)=:-K(t-\tau,k).$$
Let us consider one given fixed mode $k$ and on that account for the moment omit the $k$-dependency and introduce the notations $\phi_i(t) := \widehat{\rho}^{h_i}(k,t), \phi_e(t) := \widehat{\rho}^{h_e}(k,t)$ and $a_i(t) := \widetilde h_{i,in}(k,kt), a_e(t) := \widetilde h_{e,in}(k,kt)$. Then system \eqref{den} can be rewritten as
\begin{align}
\begin{split}
\phi_e(t)&=a_e(t)- \int_0^t K(t-\tau)[\phi_i(\tau) - \phi_e(\tau) ] \, \mathrm{d}\tau, \\
\phi_i(t)&=a_i(t)+ \frac{m_e}{m_i}\int_0^t K(t-\tau)[\phi_i(\tau) - \phi_e(\tau) ] \, \mathrm{d}\tau.
\end{split}
\label{type}
\end{align}
We can now show the convergence rates of the densities giving the damping effect under certain assumptions. 

\begin{theorem}\label{TheoremDamping}
Let $K=K(t)$ be a kernel defined for $t\geq 0$, such that 
\begin{enumerate}
\item[(i)] $|K(t)| \leq C_0 e^{-2 \pi \lambda_0 t}$ for some constants $C_0, \lambda_0>0$. 
\item[(ii)] $|(1+\frac{m_e}{m_i})K^L(\xi) -1| \geq \kappa >0$ for $0\leq$ Re$(\xi)$ $\leq \Lambda$ where the index $L$ denotes the complex Laplace transform.
\end{enumerate}
Note, that property $(i)$ in this theorem ensures that the Laplace transform is well-defined. 
Let further $a_i=a_i(t), a_e=a_e(t)$ satisfy 
\begin{enumerate}
\item[(iii)]$|\frac{m_e}{m_i}a_e(t) + a_i(t)| \leq \alpha_+ e^{-2 \pi \lambda_+ t}$,
\item[(iv)]$|a_i(t)-a_e(t)| \leq \alpha_- e^{-2 \pi \lambda_- t}$
\end{enumerate}
 and let $\phi_i, \phi_e$ solve \eqref{type}. Then for any $\lambda' < \min(\lambda_+, \lambda_-, \Lambda, \lambda_0)$ 
\begin{align*}
|\phi_e(t)| \leq C_e  e^{-2 \pi \lambda' t}, \\ |\phi_i(t)| \leq C_i  e^{-2 \pi \lambda' t},
\end{align*}
where $C_i,C_e$ are constants depending on $\lambda_+, \lambda_-, \Lambda, \lambda_0, \kappa, C_0, \lambda'$. 
\end{theorem}

\begin{proof}
Let us write $\Phi_e(t) = e^{2 \pi \lambda' t} \phi_e(t), \Phi_i(t) = e^{2 \pi \lambda' t} \phi_i(t), A_e(t)= e^{2 \pi \lambda' t} a_e(t), A_i(t)=e^{2 \pi \lambda' t} a_i(t).$ Now similar as in the proof of Lemma 3.5 in \cite{Villani_steps}, we multiply \eqref{type} by $e^{2 \pi \lambda' t}$ and obtain
\begin{align}
\begin{split}
\Phi_e(t)&=A_e(t)- \int_0^t K(t-\tau)e^{2 \pi \lambda' (t-\tau)}[\Phi_i(\tau) - \Phi_e(\tau) ] \, \mathrm{d}\tau, \\
\Phi_i(t)&=A_i(t)+\frac{m_e}{m_i}\int_0^t K(t-\tau) e^{2 \pi \lambda' (t-\tau)}[\Phi_i(\tau) - \Phi_e(\tau) ]\, \mathrm{d}\tau.
\end{split}
\label{type2}
\end{align}
The functions $\Phi_e, \Phi_i, A_e, A_i$ and $K$ are defined only for $t \geq 0$. We extend the domain of these functions for negative times by setting them equal to zero for $t<0$. Then we take the
Fourier transform in the time variable. This leads to
\begin{align*}
\widehat{\Phi}_e(\omega) &= \hat{A}_e(\omega) - \int_{-\infty}^{\infty} e^{-2\pi i \omega t} \int_0^t K(t- \tau) e^{2 \pi \lambda' (t- \tau)} [\Phi_i(\tau)- \Phi_e(\tau)] \, \mathrm{d}\tau  \mathrm{d}t \\ &= \widehat{A}_e(\omega) - \int_{-\infty}^{\infty} \int_{-\infty}^{\infty} K(s)  e^{2 \pi \lambda' s } e^{- 2\pi i \omega (s + \tau)} [\Phi_i(\tau)- \Phi_e(\tau)] \, \mathrm{d}\tau  \mathrm{d}s \\&=  \widehat{A}_e(\omega) - K^L(\lambda' + i \omega) [\widehat{\Phi}_i(\omega) - \widehat{\Phi}_e(\omega)], \\ \widehat{\Phi}_i(\omega) &= \widehat{A}_i(\omega) +\frac{m_e}{m_i}K^L(\lambda' + i \omega) [\widehat{\Phi}_i(\omega) - \widehat{\Phi}_e(\omega)],
\end{align*}
where $K^L$ denotes the complex Laplace transform defined as $K^L(\xi) = \int_0^\infty e^{2\pi \xi^* s}K(s)\, \mathrm{d}s$.
Subtracting the first equation  from the second one gives
$$ \widehat{\Phi}_i(\omega)- \widehat{\Phi}_e(\omega) = \widehat{A}_i(\omega)- \widehat{A}_e(\omega) + \left(\frac{m_e}{m_i}+1 \right)K^L(\lambda' + i \omega) [\widehat{\Phi}_i(\omega) - \widehat{\Phi}_e(\omega)],$$
which is equivalent to 
$$\widehat{\Phi}_i(\omega)- \widehat{\Phi}_e(\omega) = \frac{\widehat{A}_i(\omega)- \widehat{A}_e(\omega)}{1- \left(\frac{m_e}{m_i}+1\right) K^L(\lambda' + i \omega)},$$
since we assumed $ \left(\frac{m_e}{m_i}+1\right) K^L(\lambda' + i \omega)\neq 1$. Actually, we assumed that $|\left(\frac{m_e}{m_i}+1\right)K^L(\xi) - 1| \geq \kappa>0$ and so we obtain
$$ || \widehat{\Phi}_i- \widehat{\Phi}_e ||_{L^2(\omega)} \leq \frac{||\widehat{A}_i- \widehat{A}_e||_{L^2(\omega)} }{\kappa}.$$
Therefore by Plancherel's identity and the decay assumption $(iv)$, we get
$$ || \Phi_i- \Phi_e ||_{L^2(t)} \leq \frac{||A_i- A_e||_{L^2(t)} }{\kappa} \leq \frac{\alpha_-}{\kappa \sqrt{4 \pi (\lambda_- - \lambda')}},$$
since the integral $\int_0^{\infty} e^{- 2 \pi \lambda_- t } e^{2 \pi \lambda' t} dt$ is equal to $\frac{1}{\sqrt{4 \pi (\lambda_- - \lambda')}}$.
We plug this into the following system which is equivalent to system \eqref{type2}
\begin{align}
\begin{split}
\frac{m_e}{m_i}\Phi_e(t) + \Phi_i(t) &=\frac{m_e}{m_i}A_e(t)+ A_i(t), \\
\Phi_i(t)- \Phi_e(t) &=A_i(t)- A_e(t)+ \left(\frac{m_e}{m_i}+1\right)\int_0^t K(t-\tau) e^{2 \pi \lambda' (t-\tau)}[\Phi_i(\tau) - \Phi_e(\tau) ]\, \mathrm{d}\tau,
\end{split}
\label{type3}
\end{align}
and obtain by using the estimates $(iii)$ and $(iv)$, H\"older inequality and estimate $(i)$
\begin{align}
\begin{split}
||\frac{m_e}{m_i}\Phi_e + \Phi_i||_{L^{\infty}(t)} &\leq ||\frac{m_e}{m_i}A_e + A_i||_{L^{\infty}(t)} \leq \alpha_+, \\
|| \Phi_i - \Phi_e||_{L^{\infty}(t)} &\leq || A_i - A_e||_{L^{\infty}(t)} + \left(\frac{m_e}{m_i}+1\right)|| (K e^{2 \pi \lambda' t}) *[\Phi_i - \Phi_e]||_{L^{\infty}(dt)} \\ 
&\leq \alpha_- +\left(\frac{m_e}{m_i}+1\right)|| K e^{2 \pi \lambda' t}||_{L^2(dt)} ||\Phi_i - \Phi_e||_{L^2(dt)} \\ & \leq \alpha_- +\left(\frac{m_e}{m_i}+1\right)\frac{C_0}{\sqrt{4 \pi(\lambda_0 - \lambda')}} \frac{\alpha_-}{\kappa \sqrt{4 \pi (\lambda_- - \lambda')}}.
\end{split}
\label{diffdamp}
\end{align}
From this we obtain the estimates 
\begin{align*}
|| \Phi_i ||_{L^{\infty}(t)} &= || \frac{1}{\left(\frac{m_e}{m_i}+1\right)} (\frac{m_e}{m_i}\Phi_e + \Phi_i) + \frac{1}{\left(\frac{m_e}{m_i}+1\right)}\frac{m_e}{m_i}(\Phi_i - \Phi_e)||_{L^{\infty}(t)} \\ &\leq \frac{1}{\left(\frac{m_e}{m_i}+1\right)} ||\frac{m_e}{m_i}\Phi_e + \Phi_i ||_{L^{\infty}(t)} + \frac{1}{\left(\frac{m_e}{m_i}+1\right)} \frac{m_e}{m_i}| \Phi_i - \Phi_e ||_{L^{\infty}(t)} \\ & \leq \frac{1}{\left(\frac{m_e}{m_i}+1\right)} \alpha_+ +\frac{1}{\left(\frac{m_e}{m_i}+1\right)}\frac{m_e}{m_i}\left[\alpha_- +\left(\frac{m_e}{m_i}+1\right)\frac{C_0}{\sqrt{4 \pi(\lambda_0 - \lambda')}} \frac{\alpha_-}{\kappa \sqrt{4 \pi (\lambda_- - \lambda')}} \right],
\end{align*}
\begin{align*}
|| \Phi_e ||_{L^{\infty}(t)} &= || \frac{1}{\left(\frac{m_e}{m_i}+1\right)} \left(\frac{m_e}{m_i}\Phi_e + \Phi_i \right) - \frac{1}{\left(\frac{m_e}{m_i}+1\right)} (\Phi_i - \Phi_e)||_{L^{\infty}(t)} \\  & \leq \frac{1}{\left(\frac{m_e}{m_i}+1\right)} \alpha_+ +\frac{1}{\left(\frac{m_e}{m_i}+1\right)} \left[\alpha_- +\left(\frac{m_e}{m_i}+1\right)\frac{C_0}{\sqrt{4 \pi(\lambda_0 - \lambda')}} \frac{\alpha_-}{\kappa \sqrt{4 \pi (\lambda_- - \lambda')}} \right].
\end{align*}
We also have to check if the Fourier transforms of $\Phi_e, \Phi_i$ exist and if they are in $L^2(\mathbb{R})$. This is similar to the proof of Lemma 3.5 in \cite{Villani_steps}, so we only sketch the idea. We replace $\lambda'$ by a parameter $\alpha$ varying from $- \varepsilon$ to $\lambda'$. By the integrability of $K$ and Gronwall's Lemma, $\left(\frac{m_e}{m_i}+1\right)\Phi_e + \Phi_i$ and $\Phi_i - \Phi_e$ are bounded as a function of $t$. Therefore $\Phi_e$ and $\Phi_i$ are bounded as a function of $t$. So $\phi_i(k,t) e^{- \varepsilon |k| t}$, $\phi_e(k,t) e^{- \varepsilon |k| t}$ are integrable for any $\varepsilon > 0$ and continuous as $\varepsilon \rightarrow 0.$ Then assumption $(ii)$ guarantees that the bounds are uniform in the strip $0 \leq \text{Re}(\xi) \leq \lambda'$.
\end{proof}

\noindent From \eqref{diffdamp} we see that the difference of the densities is damped which results in a damping of the electric field and the electric energy. This is what we wanted to show. 

\begin{remark}
For a $k$-dependent kernel $K(t,k)$ the proof works similarly.
The decay assumptions have to be replaced by a decay of the order $\mathcal{O}(e^{-2\pi \lambda' |k|t})$. For more details, see the one species proof given in section 3.3 of \cite{Villani_steps}.
\end{remark}

\subsection{The Landau-Penrose stability criterion}
\label{sec6.2.2.10}

\noindent One of the assumptions of Theorem \ref{TheoremDamping}  in order to have a damping effect is the Landau-Penrose stability criterion. That means we have to ensure that
$$ |\left(\frac{m_e}{m_i}+1\right)K^L(\xi) -1| \geq \kappa >0 \quad \text{for} \quad 0 \leq \text{Re}(\xi) \leq \Lambda.$$
In \cite{Villani_steps} it is computed for the one species case that 
$$ K^L ((\lambda + i \omega)|k|) = \widehat{W}(k) \int \frac{(f^{equ})'(v)}{v-w+ i \lambda} \, \mathrm{d}v, $$
using  integration by parts and the fact that $\int_0^{\infty} e^{- 2 i \pi |k| t v } e^{2 \pi (\lambda + i w) |k| t } |k| dv dt= \frac{1}{2 i \pi} \frac{1}{v-w + i \lambda}$, assuming $(f^{equ})'$ decays fast enough at infinity such that the integral exists.
Then \cite{Villani_steps} considers the limit $\lambda \rightarrow 0^+$ (and extends the statement later to a strip $0\leq \lambda \leq \Lambda$ by continuity arguments). In the limit one obtains
$$Z(k, \omega) := \widehat{W}(k) \left[ \int \frac{(f^{equ})'(v) - (f^{equ})'(\omega)}{v-\omega} \, \mathrm{d}v - i \pi (f^{equ})'(\omega)\right].$$ So the aim is to find a condition such that $Z$ does not approach $1$. If the imaginary part of $Z$ stays away from zero, then $Z$ does not approach the real value $1$. So $Z$ only approaches zero, if $(f^{equ})'(v)$ approaches zero in the imaginary part. But when the imaginary part will approach zero or equivalently if $\omega$ approaches a zeroth of $(f^{equ})'$, the real part has to stay away from one. This leads to the Penrose stability criterion 
$$ \forall_{\omega \in \mathbb{R},} \ (f^{equ})'(\omega) = 0~ \Rightarrow~ \widehat{W} (k) \int \frac{(f^{equ})'(v)}{v - \omega} \, \mathrm{d}v <1.$$
In our two species case there is an additional factor $\left(\frac{m_e}{m_i}+1\right)$ in front of $K^L$. So our Penrose stability criterion reformulates as
$$ \forall_{\omega \in \mathbb{R},} \ (f^{equ})'(\omega) = 0 ~\Rightarrow ~\left(\frac{m_e}{m_i}+1\right)\widehat{W} (k) \int \frac{(f^{equ})'(v)}{v - \omega} \, \mathrm{d}v <1.$$
Under this assumption the results of Theorem \ref{TheoremDamping} hold meaning we have shown linear Landau damping for a two species gas near Penrose stable equilibria. That this criteria can be fulfilled by a Coulomb potential of the form $\widehat{W}(k)=1/|k|^2$ can be argued similar to  section 3 and example 3.9 in \cite{Villani_steps}.

\section{Summary and Outlook}\label{Outlook}
In this paper we considered the linear Landau damping for the two-species Vlasov-Poisson system for ions and electrons near Penrose stable equilibria.
The result is an extension of the result on the linearised one species Vlasov-Poisson equation by Mouhout and Villani \cite{Villani, Villani_steps} where one species is treated as a background species. For future work we propose to extend the results on the Landau damping for the non-linear equations (non-linear Landau damping) in \cite{Villani, Villani_steps} to a two species model where both species are evolved in time and perform an extension of the numerical results found in \cite{Pirner_Num} for a neutral two-species gas to a mixture of electrons and ions.
\\

\noindent It is further of interest to consider a non collision-less two species plasma, e.g. by taking binary collisions modeled by a BGK operator \cite{BGK, PirnerBGK, Haack} into account and see whether there is in the limit of "few collisions" a qualitative agreement. This question is formally studied in \cite{Wood} for an equation for electrons interacting with each other and with background ions.

\begin{acknowledgement}
Marlies Pirner is funded by the Deutsche Forschungsgemeinschaft (DFG, German Research Foundation) under Germany's Excellence Strategy EXC 2044-390685587,
Mathematics Muenster:
Dynamics-Geometry-Structure and by the DFG project with grant no. PI 1501/2-1. % (grant no. $PI 1501/2-1$).
Lena Baumann acknowledges support by the Würzburg Mathematics Center for Communication and Interaction (WMCCI) as well as the Stiftung der Deutschen Wirtschaft.
\end{acknowledgement}

\end{document}